\documentclass[12pt,UTF8]{article}
\usepackage{iftex}
\ifpdftex
    \usepackage{CJKutf8}
\else
    \usepackage[scheme=plain]{ctex}
    
    \newenvironment{CJK*}[2]{}{}
\fi

\usepackage{amsmath,amsfonts,amsthm}
\usepackage{mathtools,bm,mathrsfs}
\usepackage{graphicx}
\usepackage[colorlinks=true, allcolors=blue]{hyperref}
\usepackage{thmtools}
\usepackage[capitalise]{cleveref}
\usepackage{standalone}
\usepackage{indentfirst}
\usepackage[affil-it]{authblk}
\usepackage{diagbox}
\usepackage{changes}
\usepackage{enumitem}
\usepackage{bbm}
\usepackage{orcidlink}

\newtheorem{theorem}{Theorem}[section]

\newtheorem{claim}[theorem]{Claim}

\newtheorem*{claim*}{Claim}

\newcommand{\textotherwise}{\text{otherwise}}

\newcommand{\CC}{\mathbb{C}}
\newcommand{\RR}{\mathbb{R}}

\newcommand{\QQ}{\mathbb{Q}}

\newcommand{\fA}{\mathfrak{A}}
\newcommand{\fB}{\mathfrak{B}}

\DeclareMathOperator{\dist}{dist}
\DeclareMathOperator{\sgn}{sgn}

\DeclarePairedDelimiter{\abs}{\lvert}{\rvert}
\DeclarePairedDelimiter{\card}{\lvert}{\rvert}
\DeclarePairedDelimiter{\set}{\lbrace}{\rbrace}
\DeclarePairedDelimiter{\norm}{\lVert}{\rVert}

\DeclarePairedDelimiter{\lrangle}{\langle}{\rangle}

\title{Generalized block diagonal Laplacian spectrum of graphs}

\author[a,1]{Yanrui Xu(\begin{CJK*}{UTF8}{gbsn}徐岩睿\end{CJK*})\footnote{email: 23010240@mail.ecust.edu.cn}}
\author[a,1]{Da Zhao(\begin{CJK*}{UTF8}{gbsn}赵达\end{CJK*})\footnote{email: zhaoda@ecust.edu.cn}~\orcidlink{0000-0002-9582-0778}}
\affil[a]{School of Mathematics, East China University of Science and Technology, 130 Meilong Road, Shanghai 200237, China.}

\date{}

\begin{document}

\maketitle

\begin{abstract}
    We reduce the $p^2$ block all-one matrices in the generalized block Laplacian spectrum of graphs to $p$ block all-one matrices in the generalized block diagonal Lapalcian spectrum of graphs introduced by Wang and the second author (\textit{Adv. Appl. Math.} 173B (2026)). 
    In this case the matrices are all real symmetric, and hence the spectrum is real, which does not hold for the generalized block Laplacian spectrum. 
    We also investigate the analogue by Hermitian adjacency matrix of digraphs. 
\end{abstract}



\footnotetext[1]{All authors are indexed in alphabetical order. All authors are co-first authors.}

\section{Introduction}

Haemers~\cite{vandam2003WhichGraphsAre,haemers2016almost} conjectures that almost all graphs are determined by their eigenvalues. 
In other words, almost all graphs have no cospectral mate. 
Vu~\cite{vu2021RecentProgressCombinatorial} raises a similar conjecture in random symmetric $\pm 1$ matrices. 
O'Rourke--Touri~\cite{orourke2016ConjectureGodsilConcerning} proves that almost all graphs are controllable. 
If a controllable graph has a generalized cospectral mate, it must be a generalized cospectral mate with an integer level $\ell \geq 2$. 
Wei Wang and the second author~\cite{wang2025AlmostAllGraphs} prove that almost all graphs have no cospectral mate with fixed level. 
However, the level could be large compared to the order of the graph. 
Wei Wang and the second author~\cite{wang2026GraphIsomorphismMultivariate} introduce the concept of multivariate graph spectrum, which unifies the adjacency spectrum, Laplacian spectrum, generalized spectrum (idiosyncratic spectrum in the sense of Tutte~\cite{MR538033}), etc., aiming to force the level to be $1$, namely to forbid cospectral mates. 
It is shown that the multivariate graph spectrum essentially reduces to the univariate graph spectrum with well-chosen parameters. 
The generalized block Laplacian spectrum, the generalized diagonal block Laplacian spectrum, and the generalized diagonal spectrum (cospectrality in this sense is equivalent to degree-similarity in~\cite{godsil2025DegreesimilarGraphs}) are introduced there. 
In~\cite{wang2026GraphIsomorphismMultivariate} cospectrality in the sense of the generalized block Laplacian spectrum is characterized. 
However, the generalized block Laplacian spectrum involves $p^2$ block all-one matrices and the spectrum is not real in general. 
In this paper, we show that almost identical characterization holds for the generalized block diagonal Laplacian spectrum, which involves only $p$ block all-one matrices and the spectrum is always real. 
We also investigate the analogue by Hermitian adjacency matrix of digraphs. 

A \emph{graph} $G$ is a pair $(V, E)$, where $V$ is the vertex set, and $E \subseteq \binom{V}{2}$ is the edge set. 
We denote an edge $(u, v) \in E$ by $uv$.
Note that here the graph is \emph{simple}, namely it has no loops, no multiedges, and the edges are undirected. 
The degree $\deg(v)$ of a vertex $v \in V$ is the number of edges incident to $v$.
The \emph{adjacency matrix} of a graph $G$ is a zero-one square matrix $A_G$ over the field $\CC$ given by
\begin{align}
    A_G(u, v) = 
    \begin{cases}
        1, & uv \in E, \\
        0, & \textotherwise. 
    \end{cases}
\end{align}
Let $G = (V_G, E_G)$ and $H = (V_H, E_H)$ be two graphs.
A \emph{graph isomorphism} from $G$ to $H$ is a map $f : V_G \to V_H$ such that $uv \in E_G$ if and only if $f(u)f(v) \in E_H$. 
We write $G \cong_{\text{iso}} H$ if $G$ and $H$ are isomorphic.  

A matrix is called \emph{integral} if all of its entries are integers.
Throughout this paper, we shall denote by $I$ and $J$ the identity matrix and the all-one matrix respectively.  
A square integral matrix $M$ is called \emph{unimodular} if $\det M=\pm1$. Note that for a unimodular matrix $M$, its inverse $M^{-1}$ is always  an integral matrix. 
Given a rectangular $n \times m$ integral matrix $M$, there exists a decomposition $M = U \Sigma V$, called \emph{Smith decomposition}, such that the followings hold. 
\begin{enumerate}
    \item The matrix $U$ is an $n \times n$ unimodular matrix;
    \item The matrix $V$ is an $m \times m$ unimodular matrix;
    \item The matrix $\Sigma$ is a $n \times m$ diagonal matrix such that the diagonal entries satisfy $d_1 \mid d_2 \mid \cdots \mid d_{\min(n,m)}$.
\end{enumerate}
We call $\Sigma$ the \emph{Smith normal form} of $M$ and $d_{\min(n,m)}(M)$ the \emph{last invariant factor} of $M$. 

Given a rectangular matrix $M$ over $\QQ$, the \emph{level} of $M$, denoted by $\ell(M)$, is the smallest positive integer $\ell$ such that $\ell M$ is an integral matrix. 

Let $G$ be a graph. 
The (adjacency) \emph{spectrum} of $G$ is the multiset of eigenvalues of $A_G$.
We say two graphs $G$ and $H$ are \emph{cospectral}, denoted by $G \cong_{S} H$, if their spectra are identical. 
We can consider more general spectrum. 
The \emph{complementary graph} $\overline{G}$ of $G = (V, E)$ is a graph $(V, \overline{E})$ such that $\overline{E} \coloneqq \binom{V}{2} \setminus E$. 
The \emph{generalized spectrum} of $G$ is the collection of the spectrum of $G$ as well as the spectrum of $\overline{G}$. 
We say $G$ and $H$ are \emph{generalized cospectral}, denoted by $G \cong_{GS} H$, if their generalized spectra are identical. 
A matrix is called regular if the sum of the entries in each row is $1$.
The \emph{walk matrix} of a graph $G = (V, E)$ is defined by
\begin{align}
    W(G) = \begin{bmatrix}
        e & Ae & \cdots & A^{n-1}e
    \end{bmatrix},
\end{align}
where $e$ is the all-one column vector and $n = \card{V}$.
A graph $G$ is called \emph{controllable} if $W(G)$ is of full rank. 
The following theorem characterizes generalized cospectral graphs.

\begin{theorem}[{\cite{johnson1980NoteCospectralGraphs}~\cite{wang2006SufficientConditionFamily}}]\label{thm:GS}
    Let $G$ and $H$ be two non-isomorphic graphs. 
    Let $A$ and $B$ be their adjacency matrices respectively. 
    If $G$ and $H$ are generalized cospectral, then there exists a regular orthogonal matrix $Q$ such that $Q^\top A Q = B$. 
    In particular, if $G$ is controllable, then $Q$ is unique and rational and $\ell(Q) \mid d_n(W(G))$.
\end{theorem}

Let $\bm{A} = (A_1, A_2, \ldots, A_k)$ be a $k$-tuple of $n \times n$ matrices. 
Let $\bm{s} = (s_1, s_2, \ldots, s_k)$ be a $k$-tuple of complex variables. 
We define
\begin{align}
    W_{\bm{A}}(\bm{s}) = W_{A_1, A_2, \ldots, A_k}(s_1, s_2, \ldots, s_k) = \sum_{i=1}^k s_i A_i.
\end{align}
Denote by $\phi_{\bm{A}}(\bm{s}; t) = \det(t I - W_{\bm{A}}(\bm{s}))$ the \emph{characteristic polynomial} of $W_{\bm{A}}(\bm{s})$. 
It is clear that $G$ and $H$ are cospectral if and only if $\phi_{\bm{A}}(\bm{s}; t) = \phi_{\bm{B}}(\bm{s}; t)$, where $\bm{A} = (A_G)$ and $\bm{B} = (A_H)$;
$G$ and $H$ are generalized cospectral if and only if $\phi_{\bm{A}}(\bm{s}; t) = \phi_{\bm{B}}(\bm{s}; t)$, where $\bm{A} = (A_G, J)$ and $\bm{B} = (A_H, J)$~\cite{MR538033,johnson1980NoteCospectralGraphs,wang2006SufficientConditionFamily}.
Suppose $G$ and $H$ be two graphs on $V$ sharing the same degree sequences. 
Without loss of generality, suppose $\deg(v) = \deg_{G}(v) = \deg_{H}(v)$ for all $v \in V$. 
Therefore, we have the degree decomposition of vertices $V = \bigsqcup_{i=1}^p V_i$.
Let $A$ and $B$ be the adjacency matrices of $G$ and $H$ respectively. 
Let $e_i \in \RR^V$, $i = 1,2, \ldots, p$ be zero-one vectors such that $e_i(v) = 1$ if and only if $v \in V_i$. 
Let $J_{i,j} = e_i e_j^\top$ for $i,j = 1, 2, \ldots, p$.
Let 
\begin{align*}
    D_i(u, v) = 
    \begin{cases}
        1, & u = v \in V_i, \\
        0, & \textotherwise 
    \end{cases}.
\end{align*}
We say two graphs $G$ and $H$ share the \emph{generalized diagonal Laplacian spectrum}, denoted by $G \cong_{GDLS} H$, if $\phi_{\bm{A}'}(\bm{s}, t) = \phi_{\bm{B}'}(\bm{s}, t)$, where $\bm{A}' = (A_G, D_1, D_2, \ldots, D_p)$ and $\bm{B}' = (A_H, D_1, D_2, \ldots, D_p)$. 
See also~\cite{godsil2025DegreesimilarGraphs}. 
We say two graphs $G$ and $H$ share the \emph{generalized block diagonal Laplacian spectrum}, denoted by $G \cong_{GBDLS} H$, if $\phi_{\bm{A}''}(\bm{s}, t) = \phi_{\bm{B}''}(\bm{s}, t)$, where $\bm{A}'' = (A_G, J_{1,1}, J_{2,2}, \ldots, J_{p,p})$ and $\bm{B}'' = (A_H, J_{1,1}, J_{2,2}, \ldots, J_{p,p})$.
We say two graphs $G$ and $H$ share the \emph{generalized block Laplacian spectrum}, denoted by $G \cong_{GBLS} H$, if $\phi_{\bm{A}'''}(\bm{s}, t) = \phi_{\bm{B}'''}(\bm{s}, t)$, where $\bm{A}''' = (A_G, J_{1,1}, J_{1,2}, \ldots, J_{p,p})$ and $\bm{B}''' = (A_H, J_{1,1}, J_{1,2}, \ldots, J_{p,p})$. 
In~\cite{wang2026GraphIsomorphismMultivariate}, the authors generalize~\cref{thm:GS} to generalized block Laplacian spectrum.

\begin{theorem}[{\cite[Theorem 2.1]{wang2026GraphIsomorphismMultivariate}}]\label{thm:GBLS}
    Let $A_1$ and $B_1$ be two real symmetric matrices of order $n$. 
    Let $e_i$, $i = 1,2, \ldots, p$ be zero-one vectors of length $n$ such that the positions of the ones are disjoint. 
    Let $J_{i,j} = e_i e_j^\top$ for $i,j = 1, 2, \ldots, p$ and $k = 1 + p^2$.
    Let $\bm{A} = (A_1, J_{1,1}, J_{1,2}, \ldots, J_{p,p})$ and $\bm{B} = (B_1, J_{1,1}, J_{1,2}, \ldots, J_{p,p})$. 
    Then the matrices $W_{\bm{A}}(\bm{s})$ and $W_{\bm{B}}(\bm{s})$ have identically the same characteristic polynomials if and only if they are similar via a fixed orthogonal matrix $Q$, independent of $\bm{s}$, such that $Q^\top A_1 Q = B_1$ and $Q^\top e_i = e_i$ for $i = 1,2, \ldots, p$.
    Moreover, $Q^\top \widetilde{W}_{A_1} = \widetilde{W}_{B_1}$, where
    \begin{align*}
        \widetilde{W}_{A} \coloneqq [e_1, Ae_1, \ldots, A^{n-1} e_1, e_2, Ae_2, \ldots, A^{n-1} e_2, \ldots, e_p, A e_p, \ldots, A^{n-1} e_p].
    \end{align*}
    In particular, if $\widetilde{W}_{A_1}$ is of full row rank, then $Q$ is unique and rational, and $\ell(Q) \mid d_n(\widetilde{W}_{A_1})$. 
\end{theorem}

It is clear that~\cref{thm:GBLS} applies to generalized block Laplacian spectrum. 

Our first main theorem is an almost identical characterization which discards the off-diagonal block all-one matrices $J_{i,j}$ with $i \neq j$. 

\begin{theorem}\label{thm:GBDLS}
    Let $G$ and $H$ be two undirected graphs on a vertex set $V$ of cardinality $n$. 
    Let $A_1$ and $B_1$ be adjacency matrices of $G$ and $H$ respectively. 
    Let $V_i$, $i = 1, 2, \ldots, p$ be a partition of $V$ and let $e_i$ be the characteristic vector for $V_i$. 
    Define $J_{i,i} = e_i e_i^\top$ for $i = 1, 2, \cdots, p$ and $k = 1 + p$. 
    Let $\bm{A} = (A_1, J_{1,1}, J_{2,2}, \cdots, J_{p,p})$ and $\bm{B} = (B_1, J_{1,1}, J_{2,2}, \cdots, J_{p,p})$. 
    Then the matrices $W_{\bm{A}}(\bm{s})$ and $W_{\bm{B}}(\bm{s})$ have identically the same characteristic polynomials if and only if there exists a fixed orthogonal matrix Q, independent of $\bm{s}$, such that $Q^\top A_1Q=B_1$ and $Q^\top e_i= e_i$ for $i = 1, 2, \cdots, p$. 
    Consequently, $Q^\top \widetilde{W}_{A_1} = \widetilde{W}_{B_1}$, where
    \begin{align*}
        \widetilde{W}_{A} \coloneqq [e_1, Ae_1, \ldots, A^{n-1} e_1, e_2, Ae_2, \ldots, A^{n-1} e_2, \ldots, e_p, A e_p, \ldots, A^{n-1} e_p].
    \end{align*}
    In particular, if $\widetilde{W}_{A_1}$ is of full row rank, then $Q$ is unique and rational, and $\ell(Q) \mid d_n(\widetilde{W}_{A_1})$.
\end{theorem}

It is clear that~\cref{thm:GBDLS} applies to generalized block diagonal Laplacian spectrum. 

We also investigate the analogue by the Hermitian adjacency matrix of digraphs. 

A digraph $G$ is a pair $(V, E)$, where $V$ is the vertex set, and $E \subseteq \set{(u, v) \in V \mid  u \neq v}$ is the directed edge set. 
The Hermitian adjacency matrix of a digraph $G$ is a Hermitian matrix $A_G$ over the field $\CC$ given by
\begin{align}
    A_G(u, v) = 
    \begin{cases}
        1, & (u, v), (v, u) \in E, \\
        i, & (u, v) \in E, (v, u) \not\in E, \\
        -i, & (u, v) \not\in E, (v, u) \in E, \\
        0, & \textotherwise. 
    \end{cases}
\end{align}

Our second main theorem is an analogue of~\cref{thm:GBLS} in digraphs. 

\begin{theorem}\label{thm:GBLS_complex}
    Let $G$ and $H$ be two directed graphs on a vertex set $V$ of cardinality $n$. 
    Let $A_1$ and $B_1$ be the Hermitian adjacency matrices of $G$ and $H$ respectively. 
    Let $V_i$, $i = 1, 2, \ldots, p$ be disjoint subsets of $V$ and let $e_i$ be the characteristic vector for $V_i$. 
    Define $J_{i,j} = e_i e_j^\top$ for $i,j = 1, 2, \cdots, p$ and $k = 1 + p^2$. 
    Let $\bm{A} = (A_1, J_{1,1}, J_{1, 2}, \ldots, J_{p, p})$ and $\bm{B} = (B_1, J_{1,1}, J_{1, 2}, \ldots, J_{p, p})$.
    Then the matrices $W_{\bm{A}}(\bm{s})$ and $W_{\bm{B}}(\bm{s})$ have identically the same characteristic polynomials if and only if they are similar via a fixed unitary matrix $Q$, independent of $\bm{s}$, such that $Q^\dagger A_1 Q = B_1$ and $Q^\dagger e_i = e_i$.
    Moreover, $Q^\dagger \widetilde{W}_{A_1} = \widetilde{W}_{B_1}$, where
    \begin{align*}
        \widetilde{W}_{A} \coloneqq [e_1, A e_1, \ldots, A^{n-1} e_1, e_2, A e_2, \ldots, A^{n-1} e_2, \ldots, e_p, A e_p, \ldots, A^{n-1} e_p].
    \end{align*}
\end{theorem}

We fail to discard the off-diagonal block all-one matrices in the case of digraphs. 

\section{Proof}


\begin{proof}[{Proof of~\cref{thm:GBDLS}}]
    On the one hand, suppose that there exists an orthogonal matrix $Q$ such that $Q^\top A_1 Q = B_1$ and $Q^\top e_i = e_i$ for $i=1,2,\cdots,p$. 
    Then $Q^\top J_{i,i} Q = Q^\top e_i e_i^\top Q = e_i e_i^\top = J_{i,i}$. 
    Hence, $Q^\top W_{\bm{A}}(\bm{s})Q = W_{\bm{B}}(\bm{s})$, and we get $\det(t I - W_{\bm{A}}(\bm{s})) = \det(t I - W_{\bm{B}}(\bm{s}))$.

    On the other hand, suppose that $\det(t I - W_{\bm{A}}(\bm{s})) = \det(t I - W_{\bm{B}}(\bm{s}))$.
    Define a quotient graph $G^* = (X, T^*)$ with the vertex set $X = \set{u_1, u_2, \ldots, u_p}$ and $u_{i} u_{\ell} \in T$ if and only if $e_i^\top A e_\ell > 0$. 
    Without loss of generality, we may assume that $G^*$ is connected. 
    Otherwise, we can form the matrix $Q$ by putting the corresponding $Q_\Lambda$ for each connected component of $G^*$ in diagonal blocks of $Q$, up to a re-indexing of vertices.
    
    Extend the orthogonal vectors $e_1, e_2, \cdots, e_p$ to a complete orthogonal basis $e_1, e_2, \cdots, e_p, e_{p+1}, \cdots, e_n$.
    Let $O = [\frac{e_1}{\norm{e_1}}, \frac{e_2}{\norm{e_2}}, \cdots, \frac{e_p}{\norm{e_p}}, \frac{e_{p+1}}{\norm{e_{p+1}}}, \cdots, \frac{e_n}{\norm{e_n}}]$ be the orthogonal matrix. 
    Then $O^\top J_{i,i} O = n_{i,i} E_{i,i}$ where $E_{i,i}$ denotes the matrix whose only nonzero entry is a 1 in the $(i,i)$-entry. 
    Hence,
    \begin{equation}\label{eq:det}
        \det(tI - x O^\top A_1 O - \sum_{i=1}^{p} s_{i,i} n_{i,i}E_{i,i}) = \det(tI - x O^\top B_1 O - \sum_{i=1}^{p} s_{i,i} n_{i,i}E_{i,i})
    \end{equation} 
    for all $x, s_{1,1}, \cdots, s_{p,p} \in \mathbb{C}$. 
    Consider the coefficient for $s_{1,1} s_{2,2} \cdots s_{p,p}$ in the full expansion of~\cref{eq:det}.
    Then we get
    \begin{align}
        (-1)^p n_{1,1} n_{2,2} \cdots n_{p,p} \det(tI - x O^\top A_1 O - \sum_{i=1}^{p} s_{i,i} n_{i,i} E_{i,i})_{[F],[F]} \\
        = (-1)^p n_{1,1} n_{2,2} \cdots n_{p,p} \det(tI - x O^\top B_1 O - \sum_{i=1}^{p} s_{i,i} n_{i,i} E_{i,i})_{[F],[F]}
    \end{align}
    where $F$ denotes the set $\{1, 2, \cdots, p\}$ and $X_{[K],[L]}$ denotes the submatrix of $X$ obtained by deleting the rows corresponding to $K$ and the columns corresponding to $L$. 
    Namely, $\det(tI - x O^\top A_1 O)_{[F],[F]} = \det(tI - x O^\top B_1 O)_{[F],[F]}$. 
    It follows that $(O^\top A_1 O)_{[F][F]}$ and $(O^\top B_1 O)_{[F][F]}$ are orthogonally similar since they are real symmetric matrices. 
    Therefore, there exists orthogonal matrices $U_A, U_B$ such that $U_A^\top (O^\top A_1 O)_{[F][F]} U_A$ and $U_B^\top (O^\top B_1 O)_{[F][F]} U_B$ are identical diagonal matrices. 
    Let $P_A =
    \begin{bmatrix}
        I_p & 0 \\
        0 & U
    \end{bmatrix}$ and $P_B =
    \begin{bmatrix}
        I_p & 0 \\
        0 & U_B
    \end{bmatrix}$. 
    Then we may assume that after similarities of the form $P$ it holds
    \begin{equation}
       \widehat{A}_1 \coloneqq P_A^\top O^\top A_1 O P_A =
       \begin{bmatrix}
            a_{1,1} & \cdots & a_{1,p} & \alpha_{p+1,1}^\top & \cdots & \alpha_{q,1}^\top \\
            \vdots & \ddots & \vdots & \vdots & \cdots & \vdots \\
            a_{p,1} & \cdots & a_{p,p} & \alpha_{p+1,p}^\top & \cdots & \alpha_{q,q}^\top \\
            \alpha_{p+1,1} & \cdots & \alpha_{p+1,p} & \lambda_{p+1} I_{m_{p+1}} & 0 & 0 \\
            \vdots & \vdots & \vdots & 0 & \ddots & 0 \\
            \alpha_{q,1} & \cdots & \alpha_{q,p} & 0 & 0 & \lambda_q I_{m_q}
        \end{bmatrix}
    \end{equation}
    and
    \begin{equation}
       \widehat{B}_1 \coloneqq P_B^\top O^\top B_1 O P_B = 
       \begin{bmatrix}
            b_{1,1} & \cdots & b_{1,p} & \beta_{p+1,1}^\top & \cdots & \beta_{q,1}^\top \\
            \vdots & \ddots & \vdots & \vdots & \cdots & \vdots \\
            b_{p,1} & \cdots & b_{p,p} & \beta_{p+1,p}^\top & \cdots & \beta_{q,q}^\top \\
            \beta_{p+1,1} & \cdots & \beta_{p+1,p} & \lambda_{p+1} I_{m_{p+1}} & 0 & 0 \\
            \vdots & \vdots & \vdots & 0 & \ddots & 0 \\
            \beta_{q,1} & \cdots & \beta_{q,p} & 0 & 0 & \lambda_q I_{m_q}
        \end{bmatrix}
    \end{equation}
    where $\lambda_{p+1},\cdots,\lambda_q$ are the distinct eigenvalues of $(\widehat{A}_1)_{[F],[F]}$ (and also of $(\widehat{B}_1)_{[F],[F]}$). 
    Note that $P^\top E_{i,i}P=E_{i,i}$ for $i=1,2,\cdots,p$ and $P \in \set{P_A, P_B}$.
    
    Consider the coefficient for $s_{2,2}\cdots s_{p,p}$ in the full expansion of~\cref{eq:det}, and we get
    \begin{align}
        \det(tI - x\widehat{A}_1 - \sum_{i=1}^{p} s_{i,i} n_{i,i} E_{i,i})_{[\{2,3,\cdots,p\}],[\{2,3,\cdots,p\}]} \\
        = \det(tI - x\widehat{B}_1 - \sum_{i=1}^{p} s_{i,i} n_{i,i} E_{i,i})_{[\{2,3,\cdots,p\}],[\{2,3,\cdots,p\}]}.
    \end{align}
    
    For $i,\ell \in F$, we define
    \begin{equation}
        \mathfrak{A}_{i,\ell} = \delta_{i,\ell}(t - s_{i,\ell} n_{i,\ell}) - x a_{i,\ell} - \sum_{j=p+1}^q \frac{x^2 \langle \alpha_{j,i}, \alpha_{j,\ell} \rangle}{t - x \lambda_j}
    \end{equation}
    and
    \begin{equation}
        \mathfrak{B}_{i,\ell} = \delta_{i,\ell}(t - s_{i,\ell} n_{i,\ell}) - x b_{i,\ell} - \sum_{j=p+1}^q \frac{x^2 \langle \beta_{j,i}, \beta_{j,\ell} \rangle}{t - x \lambda_j}
    \end{equation}
    where $\delta_{i,\ell}$ is the Kronecker delta function, namely
    \begin{equation}
        \delta_{i,\ell} = 
        \begin{cases}
            1, & i = \ell \\
            0, & i \neq \ell.
        \end{cases}
    \end{equation}
    By the elementary row transformation, we get
    \begin{equation}
        \det(tI - x\widehat{A}_1 - \sum_{i=1}^{p} s_{i,i} n_{i,i} E_{i,i})_{[\{2,3,\cdots,p\}],[\{2,3,\cdots,p\}]} = \mathfrak{A}_{1,1} \prod_{j=p+1}^q (t - x \lambda_j)^{m_j}
    \end{equation}
    and
    \begin{equation}
        \det(tI - x\widehat{B}_1 - \sum_{i=1}^{p} s_{i,i} n_{i,i} E_{i,i})_{[\{2,3,\cdots,p\}],[\{2,3,\cdots,p\}]} = \mathfrak{B}_{1,1} \prod_{j=p+1}^q (t - x \lambda_j)^{m_j}.
    \end{equation}
    Therefore, $\mathfrak{A}_{1,1} = \mathfrak{B}_{1,1}$, which implies that $a_{1,1} = b_{1,1}$ and $\langle \alpha_{j,1}, \alpha_{j,1} \rangle = \langle \beta_{j,1}, \beta_{j,1} \rangle$ for $j = p+1, p+2, \cdots, q$. 
    Similarly, we have $a_{i,i} = b_{i,i}$ and $\langle \alpha_{j,i}, \alpha_{j,i} \rangle = \langle \beta_{j,i}, \beta_{j,i} \rangle$ for $i = 1, 2, \cdots, p$ and $j = p+1, p+2, \cdots, q$.
    
    Consider the coefficient for $s_{3,3} \cdots s_{p,p}$ in the full expansion of~\cref{eq:det}.
    Then we have
    \begin{equation*}
        \det(tI - x\widehat{A}_1 - \sum_{i=1}^p s_{i,i} n_{i,i} E_{i,i})_{[\{3,\cdots,p\}][\{3,\cdots,p\}]} = (\mathfrak{A}_{1,1} \mathfrak{A}_{2,2} - \mathfrak{A}_{1,2} \mathfrak{A}_{2,1}) \prod_{j=p+1}^{q} (t - x \lambda_j)^{m_j},
    \end{equation*}
    and
    \begin{equation*}
        \det(tI - x\widehat{B}_1 - \sum_{i=1}^p s_{i,i} n_{i,i} E_{i,i})_{[\{3,\cdots,p\}][\{3,\cdots,p\}]} = (\mathfrak{B}_{1,1} \mathfrak{B}_{2,2} - \mathfrak{B}_{1,2} \mathfrak{B}_{2,1}) \prod_{j=p+1}^{q} (t - x \lambda_j)^{m_j}.
    \end{equation*}
    Therefore,
    \begin{equation}
        \mathfrak{A}_{1,1} \mathfrak{A}_{2,2} - \mathfrak{A}_{1,2} \mathfrak{A}_{2,1} = \mathfrak{B}_{1,1} \mathfrak{B}_{2,2} - \mathfrak{B}_{1,2} \mathfrak{B}_{2,1}. \label{eq:order_2}
    \end{equation}
    Consider the polynomial term of~\cref{eq:order_2}, and we have
    \begin{align*}
        & (t - s_{1,1} n_{1,1} - x a_{1,1})(t - s_{2,2} n_{2,2} - x a_{2,2}) - (-x a_{1,2})(-x a_{2,1}) \\
        = & (t - s_{1,1} n_{1,1} - x b_{1,1})(t - s_{2,2} n_{2,2} - x b_{2,2}) - (-x b_{1,2})(-x b_{2,1}).
    \end{align*}
    Hence, $a_{1,2} a_{2,1} = b_{1,2} b_{2,1}$, namely, $a_{1,2}^2 = b_{1,2}^2$. 
    Note that $a_{1,2} = \frac{e_1}{\norm{e_1}} A_1 \frac{e_2}{\norm{e_2}} \geq 0$ and $b_{1,2} = \frac{e_1}{\norm{e_1}} B_1 \frac{e_2}{\norm{e_2}} \geq 0$. 
    Therefore, $a_{1,2} = b_{1,2} \geq 0$. 
    Consider the coefficient for $\frac{1}{t - x \lambda_j}$ in~\cref{eq:order_2}, we have
    \begin{align*}
        & (t - s_{1,1} n_{1,1} - x a_{1,1})(-x^2 \langle \alpha_{j,2}, \alpha_{j,2} \rangle) + (t - s_{2,2} n_{2,2} - x a_{2,2})(-x^2 \langle \alpha_{j,1}, \alpha_{j,1} \rangle) \\
        & - (-x a_{1,2})(-x^2 \langle \alpha_{j,1}, \alpha_{j,2} \rangle) - (-x a_{2,1})(-x^2 \langle \alpha_{j,2}, \alpha_{j,1} \rangle) \\
        = & (t - s_{1,1} n_{1,1} - x b_{1,1})(-x^2 \langle \beta_{j,2}, \beta_{j,2} \rangle) + (t - s_{2,2} n_{2,2} - x b_{2,2})(-x^2 \langle \beta_{j,1}, \beta_{j,1} \rangle) \\
        & - (-x b_{1,2})(-x^2 \langle \beta_{j,1}, \beta_{j,2} \rangle) - (-x b_{2,1})(-x^2 \langle \beta_{j,2}, \beta_{j,1} \rangle)
    \end{align*}
    Hence, $a_{1,2} \langle \alpha_{j,1}, \alpha_{j,2} \rangle = b_{1,2} \langle \beta_{j,1}, \beta_{j,2} \rangle$ for $j = p+1, p+2, \cdots, q$. 
    Consider the coefficient for $\frac{1}{(t - x \lambda_j)^2}$ in~\cref{eq:order_2}, and we have
    \begin{equation*}
        x^4 (\langle \alpha_{j,1}, \alpha_{j,1} \rangle \langle \alpha_{j,2}, \alpha_{j,2} \rangle - \langle \alpha_{j,1}, \alpha_{j,2} \rangle^2) = x^4 (\langle \beta_{j,1}, \beta_{j,1} \rangle \langle \beta_{j,2}, \beta_{j,2} \rangle - \langle \beta_{j,1}, \beta_{j,2} \rangle^2)
    \end{equation*}
    Hence, $\langle \alpha_{j,1}, \alpha_{j,2} \rangle^2 = \langle \beta_{j,1}, \beta_{j,2} \rangle^2$ for $j = p+1, p+2, \cdots, q$. Therefore,
    \begin{equation}
        \begin{cases}
            a_{1,2} = b_{1,2} > 0, \\
            \langle \alpha_{j,1}, \alpha_{j,2} \rangle = \langle \beta_{j,1}, \beta_{j,2} \rangle,
        \end{cases}
        j = p+1, p+2, \cdots, q
    \end{equation}
    or
    \begin{equation}
        \begin{cases}
            a_{1,2} = -b_{1,2} = 0, \\
            \langle \alpha_{j,1}, \alpha_{j,2} \rangle = \pm \langle \beta_{j,1}, \beta_{j,2} \rangle,
        \end{cases}
        j = p+1, p+2, \cdots, q.
    \end{equation}
    Similarly, we have
    \begin{equation}\label{eq:neq0_implies_positive}
        \begin{cases}
            a_{i,\ell} = b_{i,\ell} > 0, \\
            \langle \alpha_{j,i}, \alpha_{j,\ell} \rangle = \langle \beta_{j,i}, \beta_{j,\ell} \rangle,
        \end{cases}
        j = p+1, p+2, \cdots, q
    \end{equation}
    or
    \begin{equation}
        \begin{cases}
            a_{i,\ell} = b_{i,\ell} = 0, \\
            \langle \alpha_{j,i}, \alpha_{j,\ell} \rangle = \pm \langle \beta_{j,i}, \beta_{j,\ell} \rangle,
        \end{cases}
        j = p+1, p+2, \cdots, q
    \end{equation}
    for fixed $i,\ell = 1, 2, \cdots, p$. 

    \begin{claim}\label{clm:connected_implies_equal}
        For every $1 \leq i, \ell \leq p$ and $p+1 \leq j \leq q$, we have $a_{i, \ell} = b_{i, \ell}$ and $\lrangle{\alpha_{j, i}, \alpha_{j, \ell}} = \lrangle{\beta_{j, i}, \beta_{j, \ell}}$. 
    \end{claim}
    
    \begin{proof}
        Consider the graph distance $\dist(u_i, u_\ell)$ on $G^*$ for $u_i, u_\ell \in X$. 
        We have shown that $a_{i, i} = b_{i, i}$ and $\lrangle{\alpha_{j, i}, \alpha_{j, i}} = \lrangle{\beta_{j, i}, \beta_{j, i}}$ for $1 \leq i \leq p$ and $p+1 \leq j \leq q$. 
        By~\cref{eq:neq0_implies_positive}, we have that for $u_i, u_\ell \in X$ with $\dist(u_i, u_\ell) = 1$ it holds that $a_{i, \ell} = b_{i, \ell}$ and $\lrangle{\alpha_{j, i}, \alpha_{j, \ell}} = \lrangle{\beta_{j, i}, \beta_{j, \ell}}$ with $j = p+1, p+2, \ldots, q$.
        Now suppose that the claim holds for $1 \leq i, \ell \leq p$ with $\dist(u_i, u_\ell) \leq L$. 
        Let $u_i, u_\ell$ be two vertices such that $\dist(u_i, u_\ell) = L+1$. 
        Let $u_{i_1} u_{i_2} \cdots u_{i_{L+2}}$ be a path of length $L+1$ connecting $u_i$ and $u_{\ell}$. 
        It holds that 
        \begin{equation}
            \begin{vmatrix}
                \mathfrak{A}_{i_1,i_1} & \cdots & \mathfrak{A}_{i_1,i_{L+2}} \\
                \vdots & \ddots & \vdots \\
                \mathfrak{A}_{i_{L+2},i_1} & \cdots & \mathfrak{A}_{i_t,i_{L+2}}
            \end{vmatrix} = 
            \begin{vmatrix}
                \mathfrak{B}_{i_1,i_1} & \cdots & \mathfrak{B}_{i_1,i_{L+2}} \\
                \vdots & \ddots & \vdots \\
                \mathfrak{B}_{i_{L+2},i_1} & \cdots & \mathfrak{B}_{i_t,i_{L+2}}
            \end{vmatrix}.
        \end{equation}
        Namely,
        \begin{align}\label{eq:det_identity_sigma}
            \sum_{\sigma \in S_{L+2}} \sgn(\sigma) \prod_{k=1}^{L+2} \fA_{i_k, i_{\sigma(k)}} = \sum_{\sigma \in S_{L+2}} \sgn(\sigma) \prod_{k=1}^{L+2} \fB_{i_k, i_{\sigma(k)}}.
        \end{align}
        Note that $\fA_{i_k, i_{m}} = \fB_{i_k, i_{m}}$ for $\abs{k - m} \leq L$. 
        So we can cancel the terms with $\sigma \in S_{L+2}$ such that $\sigma(1) \neq L+2$ or $\sigma(L+2) \neq 1$. 
        For the remaining terms, we consider the coefficient for $\frac{1}{t - x \lambda_j}$. 
        Then the terms containing $t$ or $s_{i,i}$ with $i = 1, 2, \ldots, p$ can be ignored. 
        Note that $a_{i_k, i_m} = 0$ for $\abs{k-m} \geq 2$. 
        In particular, $a_{i_1, i_{L+2}} = a_{i_{L+2}, i_1} = 0$. 
        So $\frac{1}{t - x \lambda_j}$ comes from $\fA_{i_1, i_{L+2}}$ or $\fA_{i_{L+2}, i_1}$ but not both. 
        Then exactly one of $\sigma(1) = L+2$ and $\sigma(L+2) = 1$ holds. 
        If $\sigma(1) = L+2$ holds, then it contributes a nonzero term (with no more fractional parts) only if $\sigma$ is the circular permutation $((L+2) (L+1) \cdots 2 1)$ and the term is $(-1)^{L+1} \prod_{k=2}^{L+2} a_{i_k, i_{k-1}} \lrangle{\alpha_{j, i_1}, \alpha_{j, i_{L+2}}}$. 
        Similarly, if $\sigma(1) = L+2$ holds, then $\sigma$ is the circular permutation $(12 \cdots (L+1) (L+2))$ and the term is $(-1)^{L+1} \prod_{k=1}^{L+1} a_{i_k, i_{k+1}} \lrangle{\alpha_{j, i_{L+2}}, \alpha_{j, i_1}}$. 
        The two terms are in fact identical by symmetry. 
        Since $a_{i_k, i_{k+1}} = b_{i_k, i_{k+1}} > 0$, we have $\lrangle{\alpha_{j, i_1}, \alpha_{j, i_{L+2}}} = \lrangle{\beta_{j, i_1}, \beta_{j, i_{L+2}}}$.  
        Note that the above holds for all $p+1 \leq j \leq q$ and $a_{i_1, i_{L+2}} = b_{i_1, i_{L+2}} = 0$. 
        We conclude that the claim holds for $1 \leq i, \ell \leq p$ with $\dist(u_i, u_\ell) \leq L+1$. 
        Since the graph $G^*$ is connected, the full claim holds by mathematical induction. 
    \end{proof}

    By~\cref{clm:connected_implies_equal} there exists an orthogonal matrix $R = 
    \begin{bmatrix}
        I_p & 0 \\
        0 & R'
    \end{bmatrix}$ such that $\widehat{A}_1 = R^\top \widehat{B}_1 R$ and $E_{i,i} = R^\top E_{i,i} R$ for $i = 1, 2, \cdots, p$. 
    Let $Q = O P_A R^\top P_B^\top O^\top$. 
    Then
    \begin{align*}
        Q^\top A_1 Q &= O P_B R  P_A^\top O^\top A_1 O P_A R^\top P_B^\top O^\top \\
        &= O P_B R \widehat{A}_1 R^\top P_B^\top O^\top = O P_B \widehat{B}_1 P_B ^\top O^\top = B_1,
    \end{align*}
    \begin{align*}
        Q^\top J_{i,i} Q &= O P_B R P_A^\top O^\top J_{i,i} O P_A R^\top P_B^\top O^\top \\
        &= O P_B R n_{i,i} E_{i,i} R^\top P_B^\top O^\top = O P_B n_{i,i} E_{i,i} P_B ^\top O^\top = J_{i,i},
    \end{align*}
    and
    \begin{align*}
        Q^\top e_i = O P_B R P_A^\top O^\top e_i = \norm{e_i} O P_B R \chi_i = \norm{e_i} O P_B R \chi_i = e_i, 
    \end{align*}
    where $\chi_i$ is the vector whose only nonzero entry is a one in the $i$-th entry. 
    Hence, $Q^\top W_{\bm{A}}(\bm{s}) Q = W_{\bm{B}}(\bm{s})$ and $Q^\top e_i = e_i$ for $i = 1, 2, \cdots, p$.
\end{proof}

\begin{proof}[{Proof of~\cref{thm:GBLS_complex}}]
    On the one hand, suppose there exists a fixed unitary matrix $Q$ such that $Q^\dagger A_1 Q = B$ and $Q^\dagger e_i = e_i$ for $i = 1, 2, \ldots, p$.
    Then $Q^\dagger J_{i, j} = J_{i, j} = J_{i, j} Q$.
    Take conjugate transpose on both sides and we have $Q^\dagger J_{i, j} Q = J_{i, j}$ for $i, j = 1, 2, \ldots, p$.
    Hence, $\det(tI - W_{\bm{A}}(\bm{s})) = \det(Q^\dagger (tI - W_{\bm{A}}(\bm{s})) Q) = \det(tI - Q^\dagger W_{\bm{A}}(\bm{s}) Q) = \det(tI - W_{\bm{B}}(\bm{s}))$.
    Meanwhile, for $m = 0, 1,\ldots, n-1$ we have
    \begin{equation}
        Q^\dagger A_1^m e_i = B_1^m Q^\dagger e_i = B_1^m e_i.
    \end{equation}
    In other words, $Q^\dagger \widetilde{W}_{A_1} = \widetilde{W}_{B_1}$.

    On the other hand, suppose that $\det(t I - W_{\bm{A}}(\bm{s})) = \det(t I - W_{\bm{B}}(\bm{s}))$.
    Extend the orthogonal vectors $e_1, e_2, \cdots, e_p$ to a complete orthogonal basis $e_1, e_2, \cdots, e_p, e_{p+1}, \cdots, e_n$.
    Let $O = [\frac{e_1}{\norm{e_1}}, \frac{e_2}{\norm{e_2}}, \cdots, \frac{e_p}{\norm{e_p}}, \frac{e_{p+1}}{\norm{e_{p+1}}}, \cdots, \frac{e_n}{\norm{e_n}}]$ be the unitary matrix. 
    Then $O^\dagger J_{i,j} O = n_{i,j} E_{i,j}$ where $E_{i,j}$ denotes the matrix whose only nonzero entry is a 1 in the $(i,j)$-entry. 
    Hence,
    \begin{equation}\label{eq:det_complex}
        \det(tI - x O^\dagger A_1 O - \sum_{i,j=1}^{p} s_{i,j} n_{i,j}E_{i,j}) = \det(tI - x O^\dagger B_1 O - \sum_{i,j=1}^{p} s_{i,j} n_{i,j}E_{i,j})
    \end{equation} 
    for all $x, s_{1,1}, \cdots, s_{p,p} \in \mathbb{C}$. 
    Consider the coefficient for $s_{1,1} s_{2,2} \cdots s_{p,p}$ in the full expansion of~\cref{eq:det_complex}.
    Then we get
    \begin{align}
        (-1)^p n_{1,1} n_{2,2} \cdots n_{p,p} \det(tI - x O^\dagger A_1 O - \sum_{i,j=1}^{p} s_{i,j} n_{i,j} E_{i,j})_{[F],[F]} \\
        = (-1)^p n_{1,1} n_{2,2} \cdots n_{p,p} \det(tI - x O^\dagger B_1 O - \sum_{i,j=1}^{p} s_{i,j} n_{i,j} E_{i,j})_{[F],[F]}
    \end{align}
    where $F$ denotes the set $\{1, 2, \cdots, p\}$ and $X_{[K],[L]}$ denotes the submatrix of $X$ obtained by deleting the rows corresponding to $K$ and the columns corresponding to $L$. 
    Namely, $\det(tI - x O^\dagger A_1 O)_{[F],[F]} = \det(tI - x O^\dagger B_1 O)_{[F],[F]}$. 
    It follows that $(O^\dagger A_1 O)_{[F][F]}$ and $(O^\dagger B_1 O)_{[F][F]}$ are unitarily similar since they are hermitian matrices. 
    Therefore, there exists unitary matrices $U_A, U_B$ such that $U_A^\dagger (O^\dagger A_1 O)_{[F][F]} U_A$ and $U_B^\dagger (O^\dagger B_1 O)_{[F][F]} U_B$ are identical diagonal matrices. 
    Let $P_A =
    \begin{bmatrix}
        I_p & 0 \\
        0 & U_A
    \end{bmatrix}$ and $P_B =
    \begin{bmatrix}
        I_p & 0 \\
        0 & U_B
    \end{bmatrix}$. 
    Then we may assume that after similarities of the form $P$ it holds
    \begin{equation}
       \widehat{A}_1 \coloneqq P_A^\dagger O^\dagger A_1 O P_A =
       \begin{bmatrix}
            a_{1,1} & \cdots & a_{1,p} & \alpha_{p+1,1}^\dagger & \cdots & \alpha_{q,1}^\dagger \\
            \vdots & \ddots & \vdots & \vdots & \cdots & \vdots \\
            a_{p,1} & \cdots & a_{p,p} & \alpha_{p+1,p}^\dagger & \cdots & \alpha_{q,q}^\dagger \\
            \alpha_{p+1,1} & \cdots & \alpha_{p+1,p} & \lambda_{p+1} I_{m_{p+1}} & 0 & 0 \\
            \vdots & \vdots & \vdots & 0 & \ddots & 0 \\
            \alpha_{q,1} & \cdots & \alpha_{q,p} & 0 & 0 & \lambda_q I_{m_q}
        \end{bmatrix}
    \end{equation}
    and
    \begin{equation}
       \widehat{B}_1 \coloneqq P_B^\dagger O^\dagger B_1 O P_B = 
       \begin{bmatrix}
            b_{1,1} & \cdots & b_{1,p} & \beta_{p+1,1}^\dagger & \cdots & \beta_{q,1}^\dagger \\
            \vdots & \ddots & \vdots & \vdots & \cdots & \vdots \\
            b_{p,1} & \cdots & b_{p,p} & \beta_{p+1,p}^\dagger & \cdots & \beta_{q,q}^\dagger \\
            \beta_{p+1,1} & \cdots & \beta_{p+1,p} & \lambda_{p+1} I_{m_{p+1}} & 0 & 0 \\
            \vdots & \vdots & \vdots & 0 & \ddots & 0 \\
            \beta_{q,1} & \cdots & \beta_{q,p} & 0 & 0 & \lambda_q I_{m_q}
        \end{bmatrix},
    \end{equation}
    where $\lambda_{p+1},\cdots,\lambda_q$ are the distinct eigenvalues of $(\widehat{A}_1)_{[F],[F]}$ (and also of $(\widehat{B}_1)_{[F],[F]}$). 
    Note that $P^\dagger E_{i,j}P=E_{i,j}$ for $i,j = 1,2,\cdots,p$ and $P \in \set{P_A, P_B}$.

    Consider the coefficient for $s_{2,2}\cdots s_{p,p}$ in the full expansion of~\cref{eq:det_complex}, and we get
    \begin{align}
        \det(tI - x\widehat{A}_1 - \sum_{i,j=1}^{p} s_{i,j} n_{i,j} E_{i,j})_{[\{2,3,\cdots,p\}],[\{2,3,\cdots,p\}]} \\
        = \det(tI - x\widehat{B}_1 - \sum_{i,j=1}^{p} s_{i,j} n_{i,j} E_{i,j})_{[\{2,3,\cdots,p\}],[\{2,3,\cdots,p\}]}.
    \end{align}
    
    For $i,\ell \in F$, we define
    \begin{equation}
        \mathfrak{A}_{i,\ell} = \delta_{i,\ell}t - s_{i,\ell} n_{i,\ell} - x a_{i,\ell} - \sum_{j=p+1}^q \frac{x^2 \langle \alpha_{j,\ell}, \alpha_{j,i} \rangle}{t - x \lambda_j}
    \end{equation}
    and
    \begin{equation}
        \mathfrak{B}_{i,\ell} = \delta_{i,\ell}t - s_{i,\ell} n_{i,\ell} - x b_{i,\ell} - \sum_{j=p+1}^q \frac{x^2 \langle \beta_{j,\ell}, \beta_{j,i} \rangle}{t - x \lambda_j}
    \end{equation}
    where $\delta_{i,\ell}$ is the Kronecker delta function, namely
    \begin{equation}
        \delta_{i,\ell} = 
        \begin{cases}
            1, & i = \ell \\
            0, & i \neq \ell.
        \end{cases}
    \end{equation}
    By the elementary row transformation, we get
    \begin{equation}
        \det(tI - x\widehat{A}_1 - \sum_{i,j=1}^{p} s_{i,j} n_{i,j} E_{i,j})_{[\{2,3,\cdots,p\}],[\{2,3,\cdots,p\}]} = \mathfrak{A}_{1,1} \prod_{j=p+1}^q (t - x \lambda_j)^{m_j}
    \end{equation}
    and
    \begin{equation}
        \det(tI - x\widehat{B}_1 - \sum_{i,j=1}^{p} s_{i,j} n_{i,j} E_{i,j})_{[\{2,3,\cdots,p\}],[\{2,3,\cdots,p\}]} = \mathfrak{B}_{1,1} \prod_{j=p+1}^q (t - x \lambda_j)^{m_j}.
    \end{equation}
    Therefore, $\mathfrak{A}_{1,1} = \mathfrak{B}_{1,1}$, which implies that $a_{1,1} = b_{1,1}$ and $\langle \alpha_{j,1}, \alpha_{j,1} \rangle = \langle \beta_{j,1}, \beta_{j,1} \rangle$ for $j = p+1, p+2, \cdots, q$. 
    Similarly, we have $a_{i,i} = b_{i,i}$ and $\langle \alpha_{j,i}, \alpha_{j,i} \rangle = \langle \beta_{j,i}, \beta_{j,i} \rangle$ for $i = 1, 2, \cdots, p$ and $j = p+1, p+2, \cdots, q$.

    Consider the coefficient for $s_{3,3} \cdots s_{p,p}$ in the full expansion of~\cref{eq:det_complex}.
    Then we have
    \begin{equation*}
        \det(tI - x\widehat{A}_1 - \sum_{i=1}^p s_{i,i} n_{i,i} E_{i,i})_{[\{3,\cdots,p\}][\{3,\cdots,p\}]} = (\mathfrak{A}_{1,1} \mathfrak{A}_{2,2} - \mathfrak{A}_{1,2} \mathfrak{A}_{2,1}) \prod_{j=p+1}^{q} (t - x \lambda_j)^{m_j},
    \end{equation*}
    and
    \begin{equation*}
        \det(tI - x\widehat{B}_1 - \sum_{i=1}^p s_{i,i} n_{i,i} E_{i,i})_{[\{3,\cdots,p\}][\{3,\cdots,p\}]} = (\mathfrak{B}_{1,1} \mathfrak{B}_{2,2} - \mathfrak{B}_{1,2} \mathfrak{B}_{2,1}) \prod_{j=p+1}^{q} (t - x \lambda_j)^{m_j}.
    \end{equation*}
    Therefore,
    \begin{equation}
        \mathfrak{A}_{1,1} \mathfrak{A}_{2,2} - \mathfrak{A}_{1,2} \mathfrak{A}_{2,1} = \mathfrak{B}_{1,1} \mathfrak{B}_{2,2} - \mathfrak{B}_{1,2} \mathfrak{B}_{2,1}. \label{eq:order_2_complex}
    \end{equation}
    Consider the polynomial term of~\cref{eq:order_2_complex}, and we have
    \begin{align*}
        & (t - s_{1,1} n_{1,1} - x a_{1,1})(t - s_{2,2} n_{2,2} - x a_{2,2}) - (-s_{1,2} n_{1,2} - x a_{1,2})(-s_{2,1} n_{2,1} - x a_{2,1}) \\
        = & (t - s_{1,1} n_{1,1} - x b_{1,1})(t - s_{2,2} n_{2,2} - x b_{2,2}) - (-s_{1,2} n_{1,2} - x b_{1,2})(-s_{2,1} n_{2,1} - x b_{2,1}).
    \end{align*}
    Hence, $a_{1,2} = b_{1,2}$ and $a_{2,1} = b_{2,1}$. 
    Similarly, we have $a_{i,\ell} = b_{i,\ell}$ for $i,\ell = 1,2,\ldots,p$.
    Consider the coefficient for $\frac{1}{t - x \lambda_j}$ in~\cref{eq:order_2_complex}, we have
    \begin{align*}
        & (t - s_{1,1} n_{1,1} - x a_{1,1})(-x^2 \langle \alpha_{j,2}, \alpha_{j,2} \rangle) + (t - s_{2,2} n_{2,2} - x a_{2,2})(-x^2 \langle \alpha_{j,1}, \alpha_{j,1} \rangle) \\
        & - (-s_{1,2} n_{1,2} - x a_{1,2})(-x^2 \langle \alpha_{j,2}, \alpha_{j,1} \rangle) - (-s_{2,1} n_{2,1} - x a_{2,1})(-x^2 \langle \alpha_{j,1}, \alpha_{j,2} \rangle) \\
        = & (t - s_{1,1} n_{1,1} - x b_{1,1})(-x^2 \langle \beta_{j,2}, \beta_{j,2} \rangle) + (t - s_{2,2} n_{2,2} - x b_{2,2})(-x^2 \langle \beta_{j,1}, \beta_{j,1} \rangle) \\
        & - (-s_{1,2} n_{1,2} - x b_{1,2})(-x^2 \langle \beta_{j,2}, \beta_{j,1} \rangle) - (-s_{2,1} n_{2,1} - x b_{2,1})(-x^2 \langle \beta_{j,1}, \beta_{j,2} \rangle)
    \end{align*}
    Hence, $\langle \alpha_{j,1}, \alpha_{j,2} \rangle = \langle \beta_{j,1}, \beta_{j,2} \rangle$ and $\langle \alpha_{j,2}, \alpha_{j,1} \rangle = \langle \beta_{j,2}, \beta_{j,1} \rangle$ for $j = p+1, p+2, \cdots, q$. 
    Similarly, we have $\langle \alpha_{j,i}, \alpha_{j,\ell} \rangle = \langle \beta_{j,i}, \beta_{j,\ell} \rangle$ for $j = p+1,p+2,\ldots,q$ and $i,\ell = 1,2,\ldots,p$.
    In other words, the inner products among the $p$ vectors $\alpha_{j,1}, \alpha_{j,2}, \ldots, \alpha_{j,p}$ are identical to those among the $p$ vectors $\beta_{j,1}, \beta_{j,2}, \ldots, \beta_{j,p}$ for fixed $j = p+1,p+2,\ldots,q$.
    So there exists an unitary transformation $V_j$ such that $V_j^\dagger \beta_{j,i} = \alpha_{j,i}$ for $i = 1,2,\ldots,p$.
    Consider the unitary matrix
    \begin{equation}
        R = 
        \begin{bmatrix}
            I_p & 0 & \cdots & 0 \\
            0 & V_{p+1} & 0 & 0 \\
            0 & 0 & \ddots & 0 \\
            0 & 0 & 0 & V_q
        \end{bmatrix}.
    \end{equation}
    Then we have $\hat{A}_1 = R^\dagger \hat{B}_1 R$ and $R^\dagger E_{i,j} R = E_{i,j}$ for $i,j = 1,2,\ldots,p$.
    Let $Q = O P_A R^\dagger P_B^\dagger O^\dagger$.
    Then
    \begin{align*}
        Q^\dagger A_1 Q &= O P_B R P_A^\dagger O^\dagger A_1 O P_A R^\dagger P_B^\dagger O^\dagger \\
        &= O P_B R \widehat{A}_1 R^\dagger P_B^\dagger O^\dagger = O P_B \widehat{B}_1 P_B ^\dagger O^\dagger = B_1,
    \end{align*}
    \begin{align*}
        Q^\dagger J_{i,j} Q &= O P_B R P_A^\dagger O^\dagger J_{i,j} O P_A R^\dagger P_B^\dagger O^\dagger \\
        &= O P_B R n_{i,j} E_{i,j} R^\dagger P_B^\dagger O^\dagger = O P_B n_{i,j} E_{i,j} P_B ^\dagger O^\dagger = J_{i,j},
    \end{align*}
    and
    \begin{align*}
        Q^\dagger e_i = O P_B R P_A^\dagger O^\dagger e_i = \norm{e_i} O P_B R P_A^\dagger \chi_i = \norm{e_i} O P_B R \chi_i = e_i, 
    \end{align*}
    where $\chi_i$ is the vector whose only nonzero entry is a one in the $i$-th entry. 
    Hence, $Q^\dagger W_{\bm{A}}(\bm{s}) Q = W_{\bm{B}}(\bm{s})$ and $Q^\dagger e_i = e_i$ for $i = 1, 2, \cdots, p$.
\end{proof}

\section*{Acknowledgements}

Da ZHAO is supported in part by the National Natural Science Foundation of China (No. 12471324, No. 12501459, No. 12571353), and the Natural Science Foundation of Shanghai, Shanghai Sailing Program (No. 24YF2709000).  

\bibliographystyle{alpha}
\bibliography{ref}

\end{document}